\documentclass[aap,MSNbibl,citesort,seceqn,dvips]{arximspdf}

\doi{10.1214/11-AAP794}
\volume{22}
\issue{3}
\pubyear{2012}
\firstpage{1215}
\lastpage{1242}

\makeatletter

\newtheorem{theorem}{Theorem}[section]
\newtheorem{lemma}[theorem]{Lemma}
\newtheorem{proposition}[theorem]{Proposition}
\newtheorem{corollary}[theorem]{Corollary}
\newproclaim{definition}[theorem]{Definition}

\newproclaim{Remarks}{Remarks}
\newproclaim{Remark}{Remark}

\newcommand{\N}{{\mathbb{N}}}
\newcommand{\Z}{{\mathbb{Z}}}
\newcommand{\T}{{\mathbb{T}}}
\newcommand{\zd}{\Z^d}
\newcommand{\R}{{\mathbb{R}}}
\newcommand{\rd}{\R^d}
\newcommand{\C}{{\mathbb{C}}}
\newcommand{\ra}{\rightarrow}
\newcommand{\LRa}{\Leftrightarrow}
\newcommand{\ov}{\overline}
\newcommand{\tl}{\widetilde}
\newcommand{\Llra}{\Longleftrightarrow}
\newcommand{\lan}{\langle}
\newcommand{\ran}{\rangle}
\newcommand{\st}{\stackrel}
\newcommand{\nab}{\nabla}
\newcommand{\sub}{\subset}
\newcommand{\bsh}{\setminus}
\newcommand{\half}{\frac{1}{2}}
\newcommand{\Gm}{\Gamma}
\newcommand{\del}{\delta}
\newcommand{\D}{\Delta}
\newcommand{\e}{\varepsilon}
\newcommand{\tht}{\theta}
\newcommand{\lm}{\lambda}
\newcommand{\s}{\sigma}
\newcommand{\vp}{\varphi}
\newcommand{\W}{\Omega}

\newcommand{\cF}{\mathcal{F}}
\newcommand{\cG}{\mathcal{G}}
\newcommand{\cI}{\mathcal{I}}
\newcommand{\cK}{\mathcal{K}}
\newcommand{\cN}{\mathcal{N}}
\newcommand{\cP}{\mathcal{P}}
\newcommand{\cV}{\mathcal{V}}

\newcommand{\dv}{\operatorname{div}}

\makeatother

\begin{document}
\begin{frontmatter}

\title{Stochastic shear thickening fluids:
Strong convergence of the Galerkin approximation and~the~energy
equality}
\runtitle{Stochastic shear thickening fluids}

\begin{aug}
\author[A]{\fnms{Nobuo} \snm{Yoshida}\corref{}\thanksref{t1}\ead[label=e1]{nobuo@math.kyoto-u.ac.jp}}
\runauthor{N. Yoshida}
\affiliation{Division of Mathematics Graduate School
of Science Kyoto University}
\address[A]{Division of Mathematics Graduate School\\
\quad of Science\\
Kyoto University\\
Kyoto 606-8502\\
Japan\\
\printead{e1}} 
\end{aug}

\thankstext{t1}{Supported in part by JSPS Grant-in-Aid for Scientific
Research, Kiban (C) 21540125.}

\received{\smonth{2} \syear{2011}}
\revised{\smonth{6} \syear{2011}}

%
\begin{abstract}
We consider a stochastic partial differential equation (SPDE) which
describes the velocity field of a viscous, incompressible non-Newtonian
fluid subject to a random force. Here, the extra stress tensor of the
fluid is given by a polynomial of degree $p-1$ of the rate of strain
tensor, while the colored noise is considered as a random force. We
focus on the shear thickening case, more precisely, on the case $p \in
[1 +{d \over2}, {2d \over d-2})$, where $d$ is the dimension of the
space. We prove that the Galerkin scheme approximates the velocity
field in a~strong sense. As a consequence, we establish the energy
equality for the velocity field.
\end{abstract}

%

\end{frontmatter}

\section{Introduction}\label{intro}

We consider a viscous, incompressible fluid
whose motion is subject to a random force.
The container
of the fluid is supposed to be the torus
$\T^d =(\R/\Z)^d \cong[0,1]^d$ as a part of idealization.
For a differentiable vector field $v\dvtx\T^d \ra\rd$, which is interpreted
as the velocity field of the fluid,
we denote the \textit{rate of strain tensor} by
%
%
\begin{equation}\label{ev}
e(v)
= \biggl( {\partial_i v_j+\partial_j v_i \over2}\biggr)\dvtx
\T^d \ra\rd\otimes\rd.
\end{equation}
We assume that
the extra stress tensor
\[
\tau(v)\dvtx\T^d \ra\rd\otimes\rd
\]
depends on $e(v)$ polynomially. More precisely,
for $\nu>0$ (the kinematic viscosity) and $p>1$,
%
%
\begin{equation}\label{tv}
\tau( v ) = 2\nu\bigl(1 +|e(v)|^2\bigr)^{(p-2)/2}e(v).
\end{equation}
The linearly dependent case $p=2$ is the \textit{Newtonian fluid},
which is described by the Navier--Stokes equation, the special case
of~(\ref{NS1}) and~(\ref{NS2}) below.
On the other hand, both
the \textit{shear thinning} ($p<2$) and
the \textit{shear thickening} ($p>2$) cases are considered in many fields in
science and engineering. For example, shear thinning fluids are used
for automobile engine oil and pipeline for crude oil transportation,
while applications of shear thickening fluids can be found in
modeling of body armor and automobile four wheel driving systems.

We now\vspace*{1pt} explain the outline of the present paper before going through
precise definitions; cf. Sections~\ref{secfunct}--\ref{secGal} below.
The velocity field of the fluid
$X_t\dvtx\T^d \ra\R^d$ at time $t>0$, given $X_0$
is described by the following SPDE:
%
%
\begin{eqnarray}
\label{NS10}
\dv X_t&=&0,
\\
\label{NS20}
\partial_tX_t+(X_t \cdot\nab) X_t &=& -\nab\Pi_t +
\dv\tau(X_t) +\partial_t W_t.
\end{eqnarray}
Here, and in what follows,
%
%
\begin{equation}\label{divtau}
u \cdot\nab=\sum_{j=1}^du_j \partial_j
\quad\mbox{and}\quad
\dv\tau(u)=\Biggl( \sum_{j=1}^d\partial_j\tau_{ij}(u) \Biggr)_{i=1}^d
\end{equation}
for $u\dvtx\T^d \ra\R^d$.
Both the velocity field
$X_t\dvtx\T^d \ra\R^d$ and the pressure field
$\Pi_t\dvtx\T^d \ra\R$ are
the unknown process in the SPDE. The Brownian motion~$W_t$ with values in $L_2 (\T^d \ra\rd)$
(the set of vector fields on $\T^d$ with $L_2$ components)
is added as the random force.
Note also that the SPDE
(\ref{NS10}) and~(\ref{NS20}) for the case
$p=2$ is the stochastic Navier--Stokes equation~\cite{Fl08,FG95}.

In~\cite{TY10}, the following results are obtained for the
SPDE~(\ref{NS10}) and~(\ref{NS20}) in consistency with
the PDE case with nonrandom force~\cite{MNRR96}.
%
\begin{itemize}
\item There exist weak solutions for $p \in I_d$, where $I_d$ is
defined as follows: by introducing\vspace*{2pt} $p_1(d)={3d \over d+2}
\vee{3d -4 \over d}$, $p_2(d)={2d \over d-2}$ and $p_3(d)={3d-8
+\sqrt{9d^2+64} \over2d}$, $I_d=(p_1(d),\infty)$ for $2 \le d \le8$,
$I_d=(p_1(d), p_2(d)) \cup(p_3(d), \infty)$ for $d=9$ and $I_d=(p_3(d),
\infty)$ for $d \ge10$.
\item
The pathwise uniqueness of the solution holds
for $p \ge1+{d \over2}$.
\end{itemize}
We refer the readers to~\cite{TY10}, Theorems 2.1.3 and 2.2.1, for more
details of the above results.

In the case of stochastic Navier--Stokes equation, that is,
the SPDE~(\ref{NS10}) and~(\ref{NS20}) with $p=2$, the 2D (two-dimensional)
case is much better understood than the higher-dimensional case.
In particular, the weak solution is unique, which
turns out to be a strong solution~\cite{Kuk06}. It is also
known that the unique solution satisfies the energy \textit{equality},
rather than merely an inequality as in the other dimensions
\cite{Fl08,Kuk06}. We note that these nice properties
of the solution are obtained via the fact that,
for the 2D stochastic Navier--Stokes equation, the Galerkin
approximation (cf. Section~\ref{secGal} below)
converges strongly enough.\vadjust{\goodbreak}

Two progresses are made in this paper.

First is the generality.
The above-mentioned nice properties possessed by
the 2D stochastic Navier--Stokes equation
are carried over to the SPDE~(\ref{NS10}) and~(\ref{NS20})
with $p \in[1 +{d \over2}, {2d \over d-2})$.
We will do so by showing that
the associated Galerkin
approximation converges strongly enough.

The second progress made in this paper
is that the method to prove the strong convergence of the
Galerkin approximation is more direct than the
ones previously used for 2D stochastic Navier--Stokes equation,
for example,~\cite{Kuk06}. Our proof is based essentially
only on the Gronwall's lemma. In particular,
we do not need any compact embedding
theorem for Sobolev-type spaces (e.g.,~\cite{Kuk06}, page 9,
Lem\-ma~2.5).

In the rest of this section, we introduce a series of definitions
which we need to state our results precisely.

\subsection{Function spaces} \label{secfunct}

Let\vspace*{1pt} $\cV$ be the set of $\rd$-valued divergence free, mean-zero
trigonometric polynomials, that is, the set of
$v\dvtx\T^d \ra\R^d$ of the following
form:
%
%
\begin{equation}\label{cV}
v (x)=\sum_{z \in\zd\bsh\{0 \}}\widehat{v}_z\psi_z (x),\qquad
x \in\T^d,
\end{equation}
where
$\psi_z (x)=\exp(2 \pi{\mathbf i}z \cdot x)$ and the coefficients
$\widehat{v}_z \in\C^d$, $z \in\zd$ satisfy
%
%
\begin{eqnarray}
\label{cz=0}
\widehat{v}_z &= & 0
\qquad\mbox{except for finitely many $z$,} \\
\label{vreal}
\ov{\widehat{v}}_z &= & \widehat{v}_{-z}
\qquad\mbox{for all $z$}, \\
\label{kcz=0}
z \cdot\widehat{v}_z &= & 0
\qquad\mbox{for all $z$}.
\end{eqnarray}
Note that~(\ref{kcz=0}) implies that
\[
\dv v =0\qquad
\mbox{for all $v \in\cV$.}
\]
For $\alpha\in\R$ and $v \in\cV$ we define
\[
(1-\D)^{\alpha/2}v = \sum_{z \in\zd}(1+4\pi^2 |z|^2)^{\alpha
/2}\widehat{v}_z\psi_z.
\]
We equip the torus $\T^d$ with the Lebesgue measure.
For $p \in[1,\infty)$ and $\alpha\in\R$, we introduce
%
%
\begin{equation}\label{Va}
V_{p,\alpha} = \mbox{the completion of $\cV$ with respect to the
norm $ \| \cdot\|_{p,\alpha}$},
\end{equation}
where
%
%
\begin{equation}\label{normVa}
\| v \|_{p,\alpha}^p=\int_{\T^d}|(1-\D)^{\alpha/2}v|^p.
\end{equation}
Then,
%
%
\begin{equation}\label{cptimb}
V_{p,\alpha+\beta} \sub V_{p,\alpha}\qquad
\mbox{for $1\le p <\infty$, $\alpha\in\R$ and $\beta>0$},
\end{equation}
and
the inclusion $V_{p, \alpha+\beta} \ra V_{p,\alpha}$ is compact
if $1 <p <\infty$ (\cite{Ta96}, page 23, (6.9)).

\subsection{The noise}
We need the following definition.
%
\begin{definition}\label{DefHBM}
Let
$\Gm\dvtx V_{2,0} \ra V_{2,0}$ be a self-adjoint, nonnegative definite
operator of trace class.
A random variable $(W_t)_{t \ge0}$ with values in $C([0,\infty) \ra V_{2,0})$
is called a \textit{$V_{2,0}$-valued Brownian motion} with the covariance
operator $\Gm$ [abbreviated by BM($V_{2,0},\Gm$) below]
if, for each $\vp\in V_{2,0}$ and $0 \le s <t$,
\[
E [ \exp({\mathbf i}\lan\vp, W_t -W_s\ran)| (W_u)_{u \le s} ]
=\exp\biggl(-{t-s \over2}\lan\vp, \Gm\vp\ran\biggr)\qquad
\mbox{a.s.}
\]
%
\end{definition}

\subsection{The SPDE}
Given an initial velocity $X_0=\xi\in V_{2,0}$, the (random) time evolution
of the velocity field $X=(X_t)_{t \ge0}$ and
the pressure field $\Pi=(\Pi_t)_{t \ge0}$ is described by
the following SPDE: for $t >0$,
%
%
\begin{eqnarray}
\label{NS1}
X_t & \in& V_{p,1} \cap V_{2,0}, \\
\label{NS2}
\partial_tX_t+(X_t \cdot\nab) X_t &=& -\nab\Pi_t +
\dv\tau(X_t) +\partial_t W_t.
\end{eqnarray}
The formal ``time derivative''
of $W_t$, a BM($V_{2,0},\Gm$), is added as the random force.
Note that~(\ref{NS1}) implies that $ \dv X_t = 0$
in the distributional sense.
As in the case of (stochastic) Navier--Stokes equation,
we will reformulate the problem~(\ref{NS1}) and~(\ref{NS2})
into the one which does not contain the pressure.
Let
%
%
\begin{equation}\label{bv}
b(v)=-(v \cdot\nab) v+ \dv\tau(v),\qquad v \in\cV.
\end{equation}
Then, by integration by parts,
%
%
\begin{equation}\label{bv2}
\lan\vp, b(v) \ran
=\lan v, (v \cdot\nab)\vp\ran- \lan e (\vp), \tau(v) \ran,\qquad
\vp\in\cV.
\end{equation}
We generalize the definition of $b(v)$ for $v \in V_{p,1} \cap V_{2,0}$ by
regarding $b(v)$ as the linear functional on $\cV$ defined
by the right-hand side of~(\ref{bv2}).
Let $\cP\dvtx L_2 (\T^d \ra\rd) \ra V_{2,0}$ be the orthogonal projection.
Then, formally,
%
%
\begin{eqnarray}\label{PNS}
\mbox{(\ref{NS2})}\quad
& \Llra&\quad
\partial_tX_t =-\nab\Pi_t +b (X_t) +\partial_t W_t \nonumber\\
& \Longrightarrow&\quad \partial_tX_t =\cP b (X_t) +\partial_t W_t
\nonumber\\[-8pt]\\[-8pt]
&&\quad\mbox{(since $X_t, W_t \in V_{2,0}$, $\cP\circ\nab\equiv
0$)}\nonumber\\
& \Llra&\quad
X_t =X_0 +\int^t_0\cP b (X_s)\,ds +W_t.\nonumber
\end{eqnarray}
We will refer to~(\ref{NS1}) and~(\ref{PNS}) as $(\mathrm{SPLF})_p$
(stochastic power law fluid). To give a more precise
definition (Definition~\ref{Defws}), we introduce a notation.
For a Banach space $S$, we will denote
by $L_{p,\mathrm{loc}} (\R_+ \ra S)$ the set of measurable functions
$u\dvtx\R_+ \ra S$ such that $\| u \|_S $ belongs to $L_p ([0,T])$
for all $T \in(0,\infty)$, with the usual identification of any two
elements which coincide a.e.
%
\begin{definition}\label{Defws}
Let $(X,W)$ be a pair of
processes such that
$W$ is a~BM$(V_{2,0},\Gm)$.
We say that $(X,W)$ is a \textit{weak solution} to $(\mathrm{SPLF})_p$ if
the following two conditions are satisfied:
\begin{longlist}
\item[(a)] Equation~(\ref{NS1}) holds in the sense that
$t \mapsto X_t$ belongs to
%
%
\begin{equation}\label{capcap}\quad
L_{p,\mathrm{loc}} (\R_+ \ra V_{p,1})
\cap L_{\infty,\mathrm{loc}} (\R_+ \ra V_{2,0})
\cap C (\R_+ \ra V_{p'\wedge2,-\beta})
\end{equation}
for $\exists\beta>0$, where $p'={p \over p-1}$.\vspace*{1pt}
\item[(b)] Equation~(\ref{PNS}) holds in the sense that
%
%
\begin{equation}\label{SNS}
\lan\vp, X_t \ran= \lan\vp, X_0 \ran
+\int^t_0 \lan\vp, b (X_s) \ran \,ds
+\lan\vp, W_t \ran
\end{equation}
for all $\vp\in\cV$ and $t \ge0$; cf.~(\ref{bv2}).
\end{longlist}
%
\end{definition}

\subsection{The Galerkin approximation} \label{secGal}
We now discuss a finite-dimensional approximation to $(\mathrm{SPLF})_p$.

For each $z \in\zd\bsh\{ 0\}$, let
$\{e_{z,j}\}^{d-1}_{j=1}$
be an orthonormal basis of the hyperplane
$\{ x \in\rd; z \cdot x =0\}$ and let
%
%
\begin{eqnarray}\label{psi}
\psi_{z,j}(x) =
\cases{
\sqrt{2}e_{z,j}\cos(2\pi z \cdot x), &\quad $j=1,\ldots,d-1$, \cr
\sqrt{2}e_{z,j -d+1}\sin(2\pi z \cdot x), &\quad
$j=d,\ldots,2d-2$,}\nonumber\\[-8pt]\\[-8pt]
&&\eqntext{x \in\T^d.}
\end{eqnarray}
Then,
\[
\bigl\{ \psi_{z,j} ; (z,j) \in(\zd\bsh\{ 0\}) \times\{ 1,\ldots,2d-2\}
\bigr\}
\]
is an orthonormal basis of $V_{2,0}$.
We also introduce
%
%
\begin{eqnarray}\label{cV^n}
\cV_n &= & \mbox{the linear span of
$\{\psi_{z,j} ; \mbox{$(z,j)$ with $z \in[-n,n]^d$}\}$},
\nonumber\\[-8pt]\\[-8pt]
\cP_n & = & \mbox{the orthogonal projection: } V_{2,0} \ra
\cV_n.\nonumber
\end{eqnarray}
Using the orthonormal basis~(\ref{psi}),
we identify $\cV_n$ with $\R^N$, $N=\dim\cV_n$.
We suppose that:

$\blacktriangleright$
$\Gm\dvtx V_{2,0} \ra V_{2,0}$ is a self-adjoint, nonnegative definite
operator of trace class such that $\D\Gm=\Gm\D$;

$\blacktriangleright$
$W=(W_t)_{t \ge0}$ be a BM$(V_{2,0},\Gm)$ defined on a probability
space $(\W, \cF, P);$ cf. Definition~\ref{Defws};

$\blacktriangleright$
$\xi$ is a $V_{2,0}$-valued random variable defined on $(\W, \cF,
P)$ such that
%
%
\begin{equation}\label{m0}
m_0 =E [ \| \xi\|_{2,0}^2]<\infty.
\end{equation}

We note that the operator $\Gm$ has the following eigenfunction expansion
[cf.~(\ref{psi})]:
%
%
\begin{equation}\label{eigenexp}
\Gm=\sum_{z,j}\gamma^{z,j}\lan\cdot, \psi_{z,j} \ran\psi_{z,j}
\qquad\mbox{with } \gamma^{z,j} =\lan\Gm\psi_{z,j}, \psi_{z,j} \ran.
\end{equation}
We also note that $\cP_n W_t$ is identified with
an $N$-dimensional Brownian motion with covariance matrix $\Gm\cP_n$.
We consider the following
approximation of~(\ref{PNS}):
%
%
\begin{equation}\label{sdeGar^n}
X_t^n=X_0^n+\int^t_0 \cP_n b (X^n_s)\,ds+\cP_n W_t,\qquad
t \ge0,
\end{equation}
where $X_0^n=\cP_n\xi$. Let
%
%
\begin{equation}\label{X^zj}
X^{n,z,j}_t=\lan X^n_t,\psi_{z,j} \ran
\end{equation}
be the $(z,j)$-coordinate of $X^n_t$. Then,~(\ref{sdeGar^n}) reads
%
%
\begin{equation}\label{sdeGar^n2}
X^{n,z,j}_t=X^{n,z,j}_0 +\int^t_0 b^{z,j} (X^n_s)\,ds
+W^{z,j}_t,
\end{equation}
where
%
%
\begin{eqnarray}\label{b^zj}
b^{z,j} (X^n_s)&=&\lan X^n_s, (X^n_s \cdot\nab)\psi_{z,j} \ran
-\lan\tau(X^n_s), e (\psi_{z,j})\ran,\nonumber\\[-8pt]\\[-8pt]
W^{z,j}_t &=& \lan W_t,\psi_{z,j} \ran.\nonumber
\end{eqnarray}
Note also that
%
%
\begin{equation}\label{znotin[-n,n]^d}
X^{n,z,j}_t \equiv0 \qquad\mbox{if $z \notin[-n,n]^d$}.
\end{equation}
Let $W_\cdot$ and $\xi$ be as above. We then define
\begin{eqnarray*}
\cG^{\xi,W}_t & = &\s( \xi, W_s, s \le t),\qquad 0 \le t <\infty,\qquad
\cG^{\xi,W}_\infty= \s\biggl( \bigcup_{t \ge0}\cG^{\xi,W}_t \biggr),
\\
\cN^{\xi,W} & = &\{ N \sub\W, \exists\tl{N} \in\cG^{\xi
,W}_\infty,
N \sub\tl{N}, P(\tl{N})=0 \}
\end{eqnarray*}
and
%
%
\begin{equation}\label{cF^Wt}
\cF^{\xi,W}_t = \s( \cG^{\xi,W}_t \cup\cN^{\xi,W} ),\qquad
0 \le t <\infty.
\end{equation}
The following existence and uniqueness result for the SDE~(\ref{sdeGar^n})
was obtained in~\cite{TY10}.
%
\begin{theorem}\label{ThmGal}
Let $W_\cdot$, $\xi$ and $\cF^{\xi,W}_t$ be as above.
Then, for each $n \ge1$,
there exists a unique process $X_\cdot$ such that:
\begin{longlist}
\item[(a)] $X^n_t$ is $\cF^{\xi,W}_t$-measurable for all $t \ge0$;
\item[(b)]~(\ref{sdeGar^n}) is satisfied;
\item[(c)] for any $T>0$,
%
%
\begin{eqnarray}\qquad
\label{ener=}
E \biggl[ \| X^n_T \|_2^2+2\int^T_0\lan e(X^n_t),\tau(X^n_t) \ran \,dt
\biggr]
&=& E [ \| X^n_0 \|_2^2 ] +\operatorname{tr}(\Gm\cP_n )T,\\
\label{apriorip}
E \biggl[ \| X^n_T \|_2^2+{1 \over C}\int^T_0\| X^n_t \|_{p,1}^p\,dt\biggr]
&\le& m_0 +\bigl(C+\operatorname{tr}(\Gm)\bigr)T <\infty,
\end{eqnarray}
where $C=C(d, p) \in(0,\infty)$.
\end{longlist}
Suppose, in addition, that $p \ge{2d \over d+2}$. Then, for any $T>0$,
%
%
\begin{equation}\label{apriori}
E \biggl[ \sup_{t \le T}\| X^n_t \|_2^2 +\int^T_0\| X^n_t
\|_{p,1}^p\,dt\biggr] \le (1+T)C' <\infty,
\end{equation}
where $C'=C'(d, p,\Gm,m_0) \in(0,\infty)$.
\end{theorem}
%

\section{The strong convergence of the Galerkin approximation
and the energy equality}
\subsection{Strong convergence of the Galerkin approximation}
We introduce
%
%
\begin{equation}\label{lm}
\lm= \cases{
0, &\quad if $d=2$, \vspace*{2pt}\cr
\displaystyle {2(3-p)^+ \over dp -3d+4}, &\quad if $d \ge3$.}
\end{equation}
All the considerations in this article will
be limited to the case $p>{3d-4 \over d}$ if $d \ge3$
so that $\lm$ makes sense.\vspace*{1pt}

For\vspace*{1pt} $p \in[ 1+{d \over2}, {2d \over d-2})$,
the solution to $(\mathrm{SPLF})_p$
is well behaved and is well approximated by the Galerkin approximation.
%
\begin{theorem}\label{ThmGalconv}
Let $\Gm$, $W$ and $\xi$ be as in Section~\ref{secGal},
and let $X^n_t$ be the unique solution to~(\ref{sdeGar^n}); cf.
Theorem~\ref{ThmGal}.
Suppose additionally that
%
%
\begin{eqnarray}
\label{<p<}
&& d=2,3,4\quad \mbox{and}\quad 1+{d \over2} \le p <{2d \over d-2};\\
&& \mbox{the operator $\Gm\D$ is of trace class;} \\
\label{xi21}
&& \mbox{the random variable $\xi$ takes values in $V_{2,1}$}\quad
\mbox{and}\nonumber\\[-8pt]\\[-8pt]
&&m_1\st{\mathrm{def}}{=}E[\| \xi\|_{2,1}^2] <\infty.\nonumber
\end{eqnarray}
Then, there exists a process $X=(X_t)_{ t\ge0}$
on $(\W, \cF, P)$ with the following properties for any $T \in
(0,\infty)$:
\begin{longlist}
\item[(a)] For any $\alpha\in[0,1)$,
$X \in C([0,\infty) \ra V_{2,\alpha})$ and
%
%
\begin{equation}\label{supZ->08}
\sup_{0 \le t \le T}\|X_t-X^n_t\|_{2,\alpha}\st{n \ra\infty
}{\longrightarrow} 0\qquad
\mbox{in probability}.
\end{equation}
\item[(b)]
Let $\alpha\in[0,1)$ if $\lm=0$ [cf.~(\ref{lm})], and
let $\alpha=1-{2\lm\over p} \in(0,1)$ if $\lm>0$.
Then, $X \in L_{2,\mathrm{loc}} ([0,\infty) \ra V_{2,1+\alpha})$ and
%
%
\begin{equation}\label{int2a->08}
\int_0^T\| X_t-X^n_t\|_{2,1+\alpha}^2\,dt
\st{n \ra\infty}{\longrightarrow} 0 \qquad\mbox{in probability}.
\end{equation}
\item[(c)] For any $\tl{p} \in[1,p)$,
$X \in L_{\tl{p},\mathrm{loc}} ([0,\infty) \ra V_{\tl{p},1})$ and
%
%
\begin{equation}\label{intp1->08}
\lim_{n \ra\infty}E \biggl[ \int_0^T\| X_t-X^n_t\|_{\tl{p},1}^{\tl
{p}}\,dt \biggr]=0.\vadjust{\goodbreak}
\end{equation}
\end{longlist}
%
\end{theorem}

We now explain the strategy for the proof of Theorem~\ref{ThmGalconv}.
Let $Z^{m,n}_t=X^m_t-X^n_t$. Then, the core of the proof is that
%
%
\begin{equation}\label{int2->08}
\sup_{0 \le t \le T}\|Z^{m,n}_t\|_2 \st{m,n \ra\infty
}{\longrightarrow} 0\qquad
\mbox{in probability}.
\end{equation}
We will prove this by a series of elementary bounds
(mainly, Gronwall's~ine\-quality) instead of
functional analytic method based on compact embedding as in~\cite{Kuk06}.
We have by It\^{o}'s formula (cf. Lemma~\ref{Lem|Zt|^2} below for
the detail),
that
%
%
\begin{eqnarray}
\label{Ito}
\| Z_t^{m,n}\|_2^2
&=&\|(\cP_m-\cP_n)\xi\|_2^2+\operatorname{tr}(\cP_m\Gm-\cP_n \Gm)t\nonumber\\
&&{} +2M^{m,n}_t
+2\int^t_0\lan Z_s^{m,n},(\cP_m-\cP_n) b (X^n_s) \ran \,ds \\
&&{} +2\int^t_0\lan Z_s^{m,n}, b (X^m_s)-b( X^n_s)\ran \,ds,\nonumber
\end{eqnarray}
where
\[
M^{m,n}_t=\int^t_0\lan(\cP_m-\cP_n)Z^{m,n}_s, dW_s \ran.
\]
On the other hand, the following bound is known
(cf. proofs of Theorem~4.29 of~\cite{MNRR96}, pages 254 and~255, and
Theorem 2.2.1 of~\cite{TY10}) for $p>{d \over2}$
there exists $C \in(0,\infty)$ such that
%
%
\begin{equation}\label{2p-d}
\lan v-w, b (v)-b (w)\ran
\le C\| \nab v \|_p^{2p/(2p-d)}\| v-w \|_2^2\qquad
\mbox{for all $v,w \in\cV$}.\hspace*{-32pt}
\end{equation}
By~(\ref{Ito}) and~(\ref{2p-d}), we observe that for $\forall t \in[0,T]$
%
%
\begin{equation}\label{Gron1}
\| Z^{m,n}_t \|_2^2 \le S^{m,n}_T
+C\int^t_0\| \nab X^m_s \|_p^{2p/(2p-d)}\| Z^{m,n}_s \|_2^2\,ds,
\end{equation}
where
\begin{eqnarray*}
S^{m,n}_T&=&
\|(\cP_m-\cP_n)\xi\|_2^2+\operatorname{tr}(\cP_m\Gm-\cP_n \Gm)T
+2\sup_{0 \le s \le T}|M^{m,n}_s | \\
&&{} +2\int^T_0|\lan Z_s,(\cP_m-\cP_n) b (X^n_s) \ran|\,ds.
\end{eqnarray*}
Thus, by Gronwall's lemma,
%
%
\begin{eqnarray}\label{Gron2}
\sup_{0 \le t \le T}\| Z^{m,n}_t \|_2^2 \le S^{m,n}_T
\exp( CR^m_T )\nonumber\\[-8pt]\\[-8pt]
&&\eqntext{\mbox{where } \displaystyle
R^m_T=\int^t_0\| \nab X^m_s \|_p^{2p /(2p-d)}\,ds.}
\end{eqnarray}
Since ${2p \over2p-d} \le p (\LRa p \ge1+{d \over2})$,
we see from~(\ref{apriori}) that
$\{R^{m}_T\}_{m \ge1}$ are tight and so are
$\{\exp( CR^{m}_T )\}_{m \ge1}$.
Therefore, the convergence~(\ref{int2->08}) follows if
%
%
\begin{equation}\label{S->0}
S^{m,n}_T \st{m,n \ra\infty}{\longrightarrow} 0
\qquad\mbox{in probability}.
\end{equation}
This is\vspace*{1pt} shown to be true for $1+{2d \over d+2} \le p <{2d \over d-2}$;
cf. Lemma~\ref{LemS^mn->0} below. It is the most
technical part of this article and requires
a series of statements and bounds. The good news here is that
each of them is elementary.
%
\begin{Remarks*}
(1) In principle, the
Galerkin approximation converges in stronger topology for larger $p$.
It is thus natural that some lower bound of
$p$ [like\vspace*{1pt} $1+{d \over2} \le p$ in~(\ref{<p<})] is required to
show a result as above. To be precise,
the bound $1+{d \over2} \le p$ is used to get~(\ref{R^mTtight}) below.
On the other hand, the upper bound on $p$
in~(\ref{<p<}), $p <{2d \over d-2}$ is assumed for a technical
reason, which unfortunately does not seem easy to get rid of.
This technical condition guarantees the continuous embedding of
$V_{2,2}$ into
$V_{p, \alpha}$ with $\alpha>1$ and assumed rather commonly
in the literature to control the $V_{p, \alpha}$-norm of
the Galerkin approximation, for example,~\cite{MNRR96}, page 222,
(3.5) and
\cite{TY10}, proof of Lemma 3.2.2. We will need $p <{2d \over d-2}$
to be able to use
(\ref{317}) below, which is shown in~\cite{TY10}.

(2)
As mentioned in the \hyperref[intro]{Introduction}, Theorem
\ref{ThmGalconv} and the following Corollary~\ref{Corwsol} can be
thought of as an extension of the well-known case of 2D stochastic
Navier--Stokes equation ($d=p=2$) (see, e.g., \cite {Kuk06}, Theorem
2.6 and its proof). The results in the direction of Theorem
\ref{ThmGalconv} and the following Corollary~\ref{Corwsol} is also
obtained in~\cite{DD02} for the 2D Navier--Stokes equation forced by
the space--time white noise. In spite of the conceptual similarity of
their result to ours, their technique, based on the Besov spaces, is
much more involved. This is for the reason that, in contrast to the
colored noise, the white noise is so rough that the solution is not
expected to be accommodated in Sobolev spaces with positive
differentiability indices.

The existence of the
weak solution to the SPDE~(\ref{NS1}) and~(\ref{NS2})
in~\cite{TY10} includes the shear thinning case ($p<2$).
However, the weak solution discussed there is \textit{not}, in general, a
function of the initial data and the Brownian motion.
On the other hand, with Theorem~\ref{ThmGalconv}, it is almost
straightforward
to construct the weak solution to $(\mathrm{SPLF})_p$
as a function of the initial data and the Brownian motion.
\end{Remarks*}
%
\begin{corollary}\label{Corwsol}
Let $\Gm$, $W$ and $\xi$ be as in Section~\ref{secGal} and suppose
additionally that~(\ref{<p<})--(\ref{xi21}) hold true. Then, the
process $X$ in Theorem~\ref{ThmGalconv}, coupled with~$W$, is a weak
solution to $(\mathrm{SPLF})_p$ such that
%
%
\begin{eqnarray}
& & \mbox{$X_0 =\xi$}; \\
& & \mbox{$X_t$ is $\cF^{\xi,W}_t$-measurable}\qquad\mbox{for all $t \ge0$}.
\end{eqnarray}
Moreover, for any $T>0$,
%
%
\begin{equation}\label{apriori2}
E \biggl[ \sup_{t \le T}\| X_t \|_2^2 +\int^T_0\| X_t \|_{p,1}^p \,dt \biggr]
\le(1+T)C <\infty,
\end{equation}
where $C=C(d,p,\Gm,m_0) <\infty$.
\end{corollary}

We will derive Corollary~\ref{Corwsol} from
(\ref{supZ->08}) and~(\ref{intp1->08}); cf. Section~\ref{pCor}.\vspace*{-2pt}

\subsection{The energy equality}
The strong convergence of the Galerkin approximation proved in
Theorem~\ref{ThmGalconv} has the following application.\vspace*{-2pt}
%
\begin{theorem}\label{Thmee}
Let $\Gm$, $W$ and $\xi$ be as in Section~\ref{secGal} and
suppose additionally that
(\ref{<p<})--(\ref{xi21}) hold true.
Then, the pathwise energy equality holds in the sense that
there exists a martingale $M$
with respect to the filtration~$\cF^{\xi, W}_t$ such that
%
%
\begin{eqnarray}\label{Mt}
\half\| X_t \|_2^2 &=&
\half\| X_0 \|_2^2 -\int^t_0\lan e(X_s),\tau(X_s) \ran
\,ds\nonumber\\[-9pt]\\[-9pt]
&&{} +\half\operatorname{tr}(\Gm) t+M_t,\qquad
t \ge0.\nonumber
\end{eqnarray}
In particular, the mean energy equality holds
%
%
\begin{eqnarray} \label{EMt}
\half E [ \| X_t \|_2^2]
&=& \half E [ \| X_0 \|_2^2]
-E \biggl[\int^t_0\lan e(X_s),\tau(X_s) \ran \,ds \biggr]\nonumber\\[-9pt]\\[-9pt]
&&{} +\half\operatorname{tr}(\Gm) t,\qquad
t \ge0. \nonumber\vspace*{-2pt}
\end{eqnarray}
%
\end{theorem}

We prove Theorem~\ref{Thmee} by~(\ref{supZ->08}),~(\ref{intp1->08}) and
(\ref{tightp}) below; cf. Section~\ref{pee}.\vspace*{-2pt}
%
\begin{Remark*}
For the 2D stochastic Navier--Stokes equation ($d=p=2$), (\ref{Mt}) and
(\ref{EMt}) become, respectively,
%
%
\begin{eqnarray}
\label{Mt2D}
\half\| X_t \|_2^2 & = &
\half\| X_0 \|_2^2 -\nu\int^t_0\| \nab X_s
\|^2_2\,ds\nonumber\\[-9pt]\\[-9pt]
&&{}+\half\operatorname{tr}(\Gm) t +M_t,\qquad
t \ge0, \nonumber \\[-2pt]
\label{EMt2D}
\half E [ \| X_t \|_2^2]
&=&\half E [ \| X_0 \|_2^2] -\nu E \biggl[\int^t_0\| \nab X_s \|^2_2 \,ds \biggr]\nonumber\\[-9pt]\\[-9pt]
&&{} +\half\operatorname{tr}(\Gm
) t,\qquad
t \ge0.\nonumber\vspace*{-2pt}
\end{eqnarray}
\end{Remark*}

\subsection{Remarks on the 2D stochastic Navier--Stokes equation}
In this subsection, we turn to the 2D stochastic Navier--Stokes
equation, that is,
the SPDE~(\ref{NS1}) and~(\ref{NS2}) for $d=p=2$.
We remark that some important results from the
literature (e.g.,~\cite{Kuk06}, Sections 2.4 and 11.1)
follow easily from the method of the present paper.

We suppose that:\vspace*{6pt}

$\blacktriangleright$
$d=p=2$;

$\blacktriangleright$
$\Gm$, $W$ and $\xi$ are as in Section~\ref{secGal};

$\blacktriangleright$
$X^n=(X^n_t)_{t \ge0}$ is the
unique solution to~(\ref{sdeGar^n}); cf. Theorem
\ref{ThmGal}.\vadjust{\goodbreak}

We also suppose that there is an $\alpha=1,2,\ldots$ such that
%
%
\begin{eqnarray}
\label{trGD^a}\qquad
& & \mbox{the operator $\Gm(-\D)^\alpha$ is of trace class;} \\
\label{xi2a}
& & \mbox{the random variable $\xi$ takes values in $V_{2,\alpha}$
and $E[\| \xi\|_{2,\alpha}^2] <\infty$.}
\end{eqnarray}
Let $X=(X_t)_{t \ge0}$ be the limit as $n \nearrow\infty$ of the
process $X^n$ as described in
Theorem~\ref{ThmGalconv}. Then, by Corollary~\ref{Corwsol}, the
pair $(X,W)$ is identified
with the unique weak solution to
the SPDE~(\ref{NS1}) and~(\ref{NS2}). Moreover, by Theorem~\ref
{Thmee}, the
process $X$ satisfies the energy equalities
(\ref{Mt2D}) and~(\ref{EMt2D}).

\begin{proposition}\label{Prop2DNS}
Under the above assumptions, it holds for any $T \in(0,\infty)$ and
$\alpha_1 <\alpha$ that
%
%
\begin{eqnarray}
\label{2Dtight}
& & \sup_{0 \le t \le T}\| X^n_t \|_{2,\alpha}^2+\int^T_0
\| X^n_t \|_{2,\alpha+1}^2\,dt,\qquad
n=1,2,\ldots \qquad\mbox{are tight};\hspace*{-35pt}\\
\label{2Dconv}
&& \sup_{0 \le t \le T}\| X^n_t -X_t
\|_{2,\alpha_1}^2
+\int^T_0 \| X^n_t -X_t\|_{2,\alpha_1+1}^2\,dt\st{n \nearrow\infty
}{\longrightarrow} 0\nonumber\hspace*{-35pt}\\[-8pt]\\[-8pt]
&&\eqntext{\qquad\mbox{in probability}.}\hspace*{-35pt}
\end{eqnarray}
Suppose, in particular, that~(\ref{trGD^a}) and~(\ref{xi2a})
are true for $\alpha=2$. Then,
the pathwise balance
relation for the enstrophy holds in the sense that
there exists a martingale $M$
with respect to the filtration $\cF^{\xi, W}_t$ such that
%
%
\begin{eqnarray}\label{Mt2}
\half\| \nab X_t \|_2^2 &=& \half\| \nab X_0 \|_2^2
-\nu\int^t_0\| \D X_s \|_2^2 \,ds\nonumber\\[-8pt]\\[-8pt]
&&{}+\half\operatorname{tr}(-\Gm\D) t+M_t,\qquad
t \ge0.\nonumber
\end{eqnarray}
As a consequence,
%
%
\begin{eqnarray} \label{EMt2}
\half E [ \| \nab X_t \|_2^2]
&=& \half E [ \| \nab X_0 \|_2^2]
-\nu E \biggl[\int^t_0 \| \D X_s \|_2^2\,ds \biggr]\nonumber\\[-8pt]\\[-8pt]
&&{}+\half\operatorname{tr}(-\Gm\D) t,\qquad
t \ge0. \nonumber
\end{eqnarray}
\end{proposition}

We will prove Proposition~\ref{Prop2DNS} by~(\ref{supZ->08}) and
(\ref{tight2}) below; cf. Section~\ref{2DNS}.

\begin{Remark*}
The mean balance relation for the enstrophy~(\ref{EMt2}) can be used
together with~(\ref{EMt2D}) to disprove Kolmogorov-type scaling law
for 2D turbulent fluids (\cite{FGHR08}, page 11, Theorem 2.9).
\end{Remark*}

\section{\texorpdfstring{Proof of Theorem \protect\ref{ThmGalconv}}{Proof of Theorem 2.1}}
Let $n,m \in\N$, $n<m$ and
%
%
\begin{equation}\label{Zt}
Z_t=Z^{m,n}_t\st{\mathrm{def}}{=}X^m_t-X^n_t.
\end{equation}
To prove Theorem~\ref{ThmGalconv}, it is enough to prove the
following properties:\vadjust{\goodbreak}
\begin{longlist}
\item[(a)]
For any $\alpha\in[0,1)$,
%
%
\begin{equation}\label{supZ->0}
\sup_{0 \le t \le T}\|Z^{m,n}_t\|_{2,\alpha} \st{m,n \longrightarrow
\infty}{\longrightarrow} 0\qquad
\mbox{in probability}.
\end{equation}
\item[(b)]
Let $\alpha\in[0,1)$ if $\lm=0$ [cf.~(\ref{lm})] and
let $\alpha=1-{2\lm\over p} \in(0,1)$ if $\lm>0$.
Then,
%
%
\begin{equation}\label{int2a->0}
\int_0^T\| Z^{m,n}_t\|_{2,1+\alpha}^2\,dt \st{m,n \longrightarrow
\infty}{\longrightarrow} 0\qquad
\mbox{in probability}.
\end{equation}
\item[(c)] For any $\tl{p} \in[1,p)$,
%
%
\begin{equation}\label{intp1->0}
E \biggl[ \int_0^T\| Z^{m,n}_t\|_{\tl{p},1}^{\tl{p}}\,dt \biggr] \st{m,n
\longrightarrow\infty}{\longrightarrow} 0.
\end{equation}
\end{longlist}
%
\subsection{\texorpdfstring{Equation (\protect\ref{supZ->0}) implies equations (\protect\ref{int2a->0}) and (\protect\ref{intp1->0})}
{Equation (3.2) implies equations (3.3) and (3.4)}}
We first prove~(\ref{int2a->0}) and~(\ref{intp1->0}) assuming~(\ref{supZ->0}).
We will also need the following fact, which can be seen
from~\cite{TY10}, proof of Lem\-ma~3.2.2.
%
\begin{lemma}\label{Lemtlp}
%
\textup{(a)}
Suppose that $p \ge2$ if $d=2$ and that $p>p_3(d)
\st{\mathrm{def}}{=}{3d-8 +\sqrt{9d^2+64} \over2d}$
if $d \ge3$ $[$note that $p_3(d) \le1+{2d \over d+2} \le1+{d \over2}]$.
Then, ${2p \over p+2\lm} >1$ [cf.~(\ref{lm})] and
%
%
\begin{equation}\label{2,2}
E \biggl[ \int^T_0\| X^n_s \|_{2,2}^{2p / (p+2\lm)} \,dt\biggr]
\le C_T <\infty.
\end{equation}

\textup{(b)}
For $2 \le p <{2d \over d-2}$ and $\tl{p} \in(1,p)$, there exists
$\alpha>1$ such that
%
%
\begin{equation}\label{317}
E \biggl[ \int^T_0\| X^n_s \|_{p,\alpha}^{\tl{p}} \,dt\biggr]
\le C_T <\infty.
\end{equation}
%
\end{lemma}
%
\begin{pf*}{Proof of~(\ref{int2a->0})}
Let $\tht={1 \over2-\alpha} \in(0,1)$.
Then, we have by interpolation that
\[
\|Z^{m,n}_t\|_{2,1+\alpha}^2 \le
\|Z^{m,n}_t\|_{2,\alpha}^{2-2\tht}\|Z^{m,n}_t\|_{2,2}^{2 \tht}
\]
and hence, that
\[
\int^T_0\|Z^{m,n}_t\|_{2,1+\alpha}^2\,dt \le S_{m,n}^{2-2\tht}I_{m,n},
\]
where
\[
S_{m,n}=\sup_{t \le T}\|Z^{m,n}_t\|_{2,\alpha} \quad\mbox{and}\quad
I_{m,n}=\int^T_0\|Z^{m,n}_t\|_{2,2}^{2 \tht} \,dt.
\]
We note\vspace*{1pt} that $2\tht\le2p \over p+2\lm$. Since
$S_{m,n}\st{m,n \longrightarrow\infty}{\longrightarrow} 0$ in
probability by~(\ref{supZ->0}) and $\{I_{m,n} \}_{m,n \ge1}$ are tight
by~(\ref{2,2}), we get~(\ref{int2a->0}).\vadjust{\goodbreak}
\end{pf*}
\begin{pf*}{Proof of~(\ref{intp1->0})}
By~(\ref{int2a->0}),
\[
\int^T_0 \| Z^{m,n}_t \|_{1,1} \,dt \st{m,n \ra\infty}\longrightarrow0
\qquad\mbox{in probability $(P)$}.
\]
Moreover, the above random variables are uniformly integrable, since
\[
E \biggl[ \biggl( \int^T_0 \| Z^{m,n}_t \|_{1,1} \,dt\biggr)^p\biggr]
\st{\mbox{{\fontsize{8.36pt}{10.36pt}\selectfont{(\ref{apriori})}}}}{\le}
C_T <\infty.
\]
Therefore:
\begin{longlist}
\item[(1)]
\begin{eqnarray*}
\\[-28pt]
\lim_{m,n \ra\infty}E \biggl[
\int^T_0 \| Z^{m,n}_t \|_{1,1} \,dt \biggr]&=&0.
\end{eqnarray*}
\end{longlist}
Let $m (\ell), n (\ell) \nearrow\infty$ be such that
\begin{longlist}
\item[(2)]
\begin{eqnarray*}
\\[-28pt]
\Phi_{\ell, t} &\st{\mathrm{def}}{=}&\bigl|Z^{m (\ell),n (\ell)}_t \bigr|
+\bigl|\nab
Z^{m (\ell),n (\ell)}_t \bigr|\\
&\st{\ell\ra\infty}{\longrightarrow}&0,\qquad
dt|_{[0,T]} \times dx \times P\mbox{-a.e.},
\end{eqnarray*}
\end{longlist}
where $dt|_{[0,T]} \times dx$ denotes the Lebesgue measure on $[0,T]
\times\T^d$. Such sequences
$m (\ell), n (\ell) $ exist by (1). The sequence $\{ \Phi_{\ell,
\cdot} \}_{\ell\ge1}$ are uniformly
integrable with respect to $dt|_{[0,T]} \times dx \times P$. In fact,
\[
E \biggl[ \int^T_0 \int_{\T^d}\Phi_{\ell, t}^p \,dt \biggr]
\st{\mbox{{\fontsize{8.36pt}{10.36pt}\selectfont{(\ref{apriori})}}}}{\le} C_T <\infty.
\]
Therefore, (2) together with this uniform integrability implies (\ref
{intp1->0}) along
the subsequence $m (\ell), n (\ell)$.
Finally, we get rid of the subsequence, since the subsequence as $m
(\ell), n (\ell)$ above can be
chosen from any subsequence of $m,n$ given in advance.
\end{pf*}

\subsection{The bound by Gronwall's lemma} \label{secGron}
We will prove~(\ref{supZ->0}) in Sections~\ref{secGron} and~\ref{secS^mn->0}.
We start with an easy It\^{o} calculus. We write
$|z|_\infty=\max_{1 \le i \le d}|z_i|$ for $z=(z_1,\ldots,z_d) \in\rd$.
%
\begin{lemma}\label{Lem|Zt|^2}
%
%
\begin{eqnarray} \label{|Zt|^2}
\| Z_t\|_2^2
&=&\|(\cP_m-\cP_n)\xi\|_2^2+\operatorname{tr}(\cP_m\Gm-\cP_n \Gm)t \nonumber\\[-2pt]
& &{} +2M^{m,n}_t
+2\int^t_0\lan Z_s,(\cP_m-\cP_n) b (X^n_s) \ran \,ds \\[-2pt]
& &{} +2\int^t_0\lan Z_s, b (X^m_s)-b( X^n_s)\ran \,ds,\nonumber
\end{eqnarray}
where
%
%
\begin{equation}\label{M^mn}
M^{m,n}_t=\mathop{\sum_{z,j}}_{n<|z|_\infty\le m}\int^t_0Z^{z,j}_s \,
dW^{z,j}_s.\vadjust{\goodbreak}
\end{equation}
%
\end{lemma}

\begin{pf}
We write
\[
Z_t=(\cP_m-\cP_n)\xi+\int^t_0\bigl(\cP_m b(X^m_s)-\cP_n b(X^n_s)\bigr)\,ds
+(\cP_m-\cP_n)W_t.
\]
Since
\[
\| Z_t\|_2^2=\sum_{z,j}|Z^{z,j}_t|^2,
\]
we compute each summand. Recall that $n <m$.
If $|z|_\infty\le n $, then
\[
Z^{z,j}_t=\int^t_0\bigl(b^{z,j}(X^m_s)-b^{z,j}(X^n_s)\bigr)\,ds,
\]
and thus,
\[
|Z^{z,j}_t|^2=2\int^t_0 Z^{z,j}_s \bigl(b^{z,j}(X^m_s)-b^{z,j}(X^n_s)\bigr)\,ds.
\]
On the other hand, if $n<|z|_\infty\le m$, then
\[
Z^{z,j}_t=\xi^{z,j}+\int^t_0b^{z,j}(X^m_s)\,ds+W^{z,j}_t.
\]
With the martingale
\[
M^{z,j}_t=\int^t_0 Z^{z,j}_s \,dW^{z,j}_s
\]
we have
\begin{eqnarray*}
|Z^{z,j}_t|^2
&=&|\xi^{z,j}_t|^2+2\int^t_0 Z^{z,j}_sb^{z,j}(X^m_s) \,ds
+2M^{z,j}_t+\gamma^{z,j}t\\
&=&|\xi^{z,j}_t|^2+2\int^t_0 Z^{z,j}_s\bigl(b^{z,j}(X^m_s)-b^{z,j}(X^n_s)\bigr)
\,ds \\
& &{} +2\int^t_0 Z^{z,j}_sb^{z,j}(X^n_s) \,ds
+2M^{z,j}_t+\gamma^{z,j}t,
\end{eqnarray*}
where $\gamma^{z,j} =\lan\Gm\psi_{z,j}, \psi_{z,j} \ran$.
Putting these together, we get
\begin{eqnarray*}
\| Z_t\|_2^2
&=&\|(\cP_m-\cP_n)\xi\|_2^2+\operatorname{tr}(\cP_m\Gm-\cP_n \Gm)t+2M^{m,n}_t
\nonumber\\
&&{}+2\mathop{\sum_{z,j}}_{n<|z|_\infty\le m}\int^t_0
Z^{z,j}_sb^{z,j}(X^n_s) \,ds \\
&&{} +2\mathop{\sum_{z,j}}_{|z|_\infty\le m}\int^t_0 Z^{z,j}_s
\bigl(b^{z,j}(X^m_s)-b^{z,j}(X^n_s)\bigr) \,ds,\nonumber
\end{eqnarray*}
which is~(\ref{|Zt|^2}).
\end{pf}
%
\begin{lemma}\label{LemGron}
Referring to Lemma~\ref{Lem|Zt|^2}, let
%
%
\begin{eqnarray}\label{S^mn}
S^{m,n}_T&=&
\|(\cP_m-\cP_n)\xi\|_2^2+\operatorname{tr}(\cP_m\Gm-\cP_n \Gm)T
+{2\sup_{0 \le s \le T}}|M^{m,n}_s |\nonumber\\[-8pt]\\[-8pt]
& &{} +2\int^T_0|\lan Z_s,(\cP_m-\cP_n) b (X^n_s) \ran|\,ds.\nonumber
\end{eqnarray}
Then, for $p > {d \over2}$,
%
%
\begin{equation}\label{Gron}
\sup_{0 \le t \le T}\|Z_t\|_2^2
\le S^{m,n}_T\exp\biggl(C \int^T_0 \|\nab X^m_s\|_p^{2p /(2p-d)}\,ds \biggr).
\end{equation}
%
\end{lemma}
\begin{pf}
The lemma follows from Lemma~\ref{Lem|Zt|^2}, the known bound
(\ref{2p-d})
and Gronwall's lemma, exactly as explained earlier; cf.~(\ref{Gron2}).
\end{pf}

\subsection{\texorpdfstring{Proof of (\protect\ref{supZ->0})}{Proof of (3.2)}} \label{secS^mn->0}
The essential part of the proof of~(\ref{supZ->0}) is the \mbox{following}.
%
\begin{lemma}\label{LemS^mn->0}
For $1+{2d \over d+2} \le p <{2d \over d-2}$,
\[
S^{m,n}_T \st{m,n \ra\infty}{\longrightarrow} 0 \qquad\mbox{in
probability},
\]
where $S^{m,n}_T$ is defined by~(\ref{S^mn}).
\end{lemma}

Most of this subsection is devoted to the proof of Lemma~\ref{LemS^mn->0}.
Using Lem\-ma~\ref{LemS^mn->0}, we will prove~(\ref{supZ->0}) at the
end of this subsection.

Referring to~(\ref{S^mn}), it is obvious that
%
%
\begin{equation}\label{xi+Gm->0}
\|(\cP_m-\cP_n)\xi\|_2^2+\operatorname{tr}(\cP_m\Gm-\cP_n \Gm)T
\longrightarrow0,\qquad
m,n \longrightarrow0.
\end{equation}
On the other hand, it is easy to prove that
%
%
\begin{equation}\label{M^mn->0}
E\Bigl[\sup_{0 \le t \le T}|M^{m,n}_t|^2\Bigr] \longrightarrow0,\qquad m,n
\longrightarrow\infty.
\end{equation}
To see this, we compute the quadratic variation of $M^{m,n}$,
\begin{eqnarray*}
\lan M^{m,n}\ran_t&=&\int^t_0 \lan(\cP_m\Gm-\cP_n\Gm
)X^m_s,X^m_s\ran \,ds \\
& \le& \| \cP_m\Gm-\cP_n\Gm\|_{2 \ra2}\int^t_0\|X^m_s \|_2^2\,ds.
\end{eqnarray*}
Here, and in what follows, we denote the
norm of the bounded operators on~$V_{p,0}$ by
%
%
\begin{equation}\label{opnorm}
\| \cdot\|_{p \ra p}.
\end{equation}
We have that
\[
\| \cP_m\Gm-\cP_n\Gm\|_{2 \ra2}^2
\le\mathop{\sum_{z,j}}_{n <|z|_\infty\le m}|\gamma^{z,j}|^2
\longrightarrow0\vadjust{\goodbreak}
\]
and that
\[
\sup_m E \biggl[ \int^t_0\|X^m_s \|_2^2\,ds\biggr] \le C_t <\infty
\]
by~(\ref{apriori}). Thus, by Doob's $L^2$-maximal inequality,
\[
E\Bigl[\sup_{0 \le t \le T}|M^{m,n}_t|^2\Bigr]
\le4 E[ \lan M^{m,n}\ran_T] \longrightarrow0.
\]
Therefore, to prove Lemma~\ref{LemS^mn->0}, it is enough to show that
%
%
\begin{equation}\label{int^T0->0}\quad
\int^T_0|\lan Z_s,(\cP_m-\cP_n) b (X^n_s) \ran|\,ds
\st{m,n \longrightarrow\infty}{\longrightarrow} 0 \qquad\mbox{in
probability},
\end{equation}
if $1+{2d \over d+2} \le p <{2d \over d-2}$.\vspace*{1pt}

The rest of this section will be devoted
to the proof of~(\ref{int^T0->0}). We start by cutting the task
into pieces. Since $(\cP_m-\cP_n)Z_s =(1-\cP_n)X^m_s$, we have
%
%
\begin{eqnarray}\label{p1p1}
|\lan Z_s,(\cP_m-\cP_n) b (X^n_s) \ran|
&=& |\lan(1-\cP_n)X^m_s, b (X^n_s) \ran| \nonumber\\[-8pt]\\[-8pt]
&\le&
\| (1-\cP_n)X^m_s\|_{p,1}\|b (X^n_s) \|_{p',-1}.\nonumber
\end{eqnarray}
With $\alpha>1$ to be specified later on, we bound the first factor of
(\ref{p1p1}) as follows:
%
%
\begin{eqnarray}\label{p1p1fac1}
\| (1-\cP_n)X^m_s\|_{p,1}
&=&\| (1-\D)^{1/2}(1-\cP_n)X^m_s\|_p \nonumber\\
&= & \bigl\| (1-\cP_n)(1-\D)^{(1-\alpha)/2}(1-\D)^{\alpha/2}X^m_s\bigr\|
_p \nonumber\\[-8pt]\\[-8pt]
& \le& \e_n \| X^m_s\|_{p,\alpha}\nonumber\\
&&\eqntext{\mbox{where }
\e_n =\bigl\|(1-\cP_n)(1-\D)^{(1-\alpha)/2}\bigr\|_{p \ra p}.}
\end{eqnarray}
As for the second factor of~(\ref{p1p1}),
we use~\cite{TY10}, (1.31) and (1.32), to get
%
%
\begin{eqnarray}\label{p1p1fac2}
\|b (X^n_s) \|_{p',-1}
& = & \| {-}(X^n_s \cdot\nab)X^n_s+ \dv\tau(X^n_s)\|_{p',-1}
\nonumber\\[-8pt]\\[-8pt]
& \le& C\|X^n_s\|_{p,1}\|X^n_s\|_2 +C(1+\|\nab
X^n_s\|_p)^{p-1}.\nonumber
\end{eqnarray}
Putting~(\ref{p1p1})--(\ref{p1p1fac2}) together, we have
\[
\int^T_0|\lan Z_s,(\cP_m-\cP_n) b (X^n_s) \ran|\,ds
\le C\e_n ( I^{m,n}_T+J^{m,n}_T),
\]
where
%
%
\begin{eqnarray}\label{IJ}
I^{m,n}_T&=&\int^T_0 \| X^m_s \|_{p,\alpha}\| X^n_s \|_{p,1}\| X^n_s
\|_2\,ds,\nonumber\\[-8pt]\\[-8pt]
J^{m,n}_T&=&\int^T_0 \| X^m_s \|_{p,\alpha}(1+\| \nab X^n_s
\|_p)^{p-1}\,ds.\nonumber
\end{eqnarray}
We will prove~(\ref{int^T0->0}) by showing that
%
%
\begin{eqnarray}\label{en->}
&\e_n \ra0 \qquad\mbox{for any $\alpha>1$};&
\\
%
%
\label{IJtight}
&\{ I^{m,n}_T\}_{m,n}, \{ J^{m,n}_T\}_{m,n}\mbox{ are tight} \qquad \mbox{for
some $\alpha>1$.}&
\end{eqnarray}
Since $(1-\D)^{(1-\alpha)/2}\dvtx V_{p,0} \ra V_{p,0}$ is compact for
any $\alpha>1$,
(\ref{en->}) follows from Lemma 3.3.2.
%
\begin{lemma}\label{LemGcpt}
Let $G\dvtx V_{p,0} \ra V_{p,0}$ be a compact operator. Then
\[
\lim_{n \ra\infty}\|(1-\cP_n)G\|_{p \ra p}=0.
\]
%
\end{lemma}
%
\begin{pf}
Since the projection $\cP_n$ corresponds to the rectangular partial summation
of the Fourier series, $\| \cP_n \|_{p \ra p}$ is bounded
in $n$ (see, e.g.,~\cite{Gra04}, page 213, Theorem 3.5.7).
Assuming this, the proof of the lemma is standard (compact uniform
convergence of a series of equi-continuous functions, which
converge on a dense set).
\end{pf}

We now turn to~(\ref{IJtight}). We will use some facts from~\cite{TY10}.
For $v \in\cV$, we introduce
%
%
\begin{eqnarray}\quad
\label{cIv}
\cI_p (v) &=& \int_{\T^d}\bigl(1+|e (v)|^2\bigr)^{(p-2)/2}|\nab e
(v)|^2,
\\
\label{cKv}
\cK(v) &= & \lan-\D v, (v \cdot\nab v) v \ran
- \lan\tau(v), e (-\D v ) \ran+\tfrac{1}{2}\operatorname{tr} (-\D\Gm\cP_n).
\end{eqnarray}
Since $|\D v | \le| \nab e (v)|$, we have
%
%
\begin{equation}\label{D<I}
\| \D v \|_2^2 \le\cI_p (v) \qquad\mbox{for $p \ge2$.}
\end{equation}

Then, we have from the proof of Lemma 3.2.3 in~\cite{TY10} that
%
%
\begin{eqnarray}
\label{K+cI}
\cK(v) +c_1 \cI_p (v) & \le&
C_1(1+\| \nab v \|^2_2)^\lm(1+\| \nab v \|_p)^p, \\
\label{EI}
E \biggl[ \int^T_0{ \cI_p (X^n_t) \over(1+\| X^n_s \|_2^2)^\lm} \,dt\biggr]
& \le& C_T <\infty.
\end{eqnarray}
%
Having prepared all the ingredients from~\cite{TY10},
our starting point to prove~(\ref{IJtight})
is the following tightness lemma (Lemma~\ref{Lemtight2}).
In fact, this tightness, together with Lemma~\ref{Lemtlp},
is enough for the proof of~(\ref{IJtight}) for $p=2;$
cf. case~1 in the proof of Lemma~\ref{LemS^mn->0} below.
%
\begin{lemma}\label{Lemtight2}
Let $p \ge1+{2d \over d+2} \ge2$. Then
%
%
\begin{equation}\label{tight2}
\sup_{0 \le t \le T}\| X^n_t\|_{2,1}
,\qquad n=1,2,\ldots,\qquad\mbox{are tight}.
\end{equation}
%
\end{lemma}
%
\begin{pf}
Note that $p \ge1+{2d \over d+2}>{3d-4 \over d}$. For $x \ge0$, let
\[
f(x)=\cases{
\displaystyle {1 \over1-\lm}(1+x)^{1-\lm}, &\quad if $\lm\neq1$, \vspace*{2pt}\cr
\ln(1+x), &\quad if $\lm= 1$.}\vadjust{\goodbreak}
\]
The condition $p \ge1+{2d \over d+2}$ guarantees that $\lm\in[0,1]$
and hence, that
\[
0 \le f(x) \ra\infty\qquad\mbox{as $x \ra\infty$}.
\]
Thus, taking~(\ref{apriori}) into account, it is enough to prove that
\begin{longlist}
\item[(1)]
\begin{eqnarray*}
\\[-28pt]
E \Bigl[\sup_{0 \le t \le T}f(\| \nab X^n_t\|_2^2) \Bigr]
\le C_T <\infty.
\end{eqnarray*}
\end{longlist}
We have by It\^o's formula that
\begin{longlist}
\item[(2)]
\begin{eqnarray*}
\\[-28pt]
f(\| \nab X^n_t \|_2^2 )
&\le& f(\| \nab X^n_0 \|_2^2 )
+N^n_t
+2\int^t_0 {\cK(X^n_s) \,ds \over(1+\| \nab X^n_s \|_2^2 )^\lm},
\end{eqnarray*}
\end{longlist}
where
\[
N^n_t = \sum_{z,j}\int^t_0
{ \D X^{n, z,j}_s \over(1+\| \nab X^n_s \|_2^2 )^\lm}\,d W^{z,j}_s;
\]
cf.~\cite{TY10}, proof of Lemma 3.2.3.
We see from~(\ref{K+cI}) that
\[
\sup_{0 \le t \le T}\int^t_0 {\cK(X^n_s) \,ds \over(1+\| \nab X^n_s \|
_2^2 )^\lm}
\le C_1\int^T_0 (1+\| \nab X^n_s \|_p)^p\,ds
\]
and hence, that
\begin{longlist}
\item[(3)]
\begin{eqnarray*}
\\[-28pt]
E \biggl[\sup_{0 \le t \le T}\int^t_0
{\cK(X^n_s) \,ds \over(1+\| \nab X^n_s \|_2^2 )^\lm} \biggr] &\le&
C_T <\infty
\end{eqnarray*}
\end{longlist}
by~(\ref{apriori}).
On the other hand, we compute
\[
\lan N^n \ran_t
=
\sum_{z,j}\int^t_0
{(\D X^{n, z,j}_s )^2 \over(1+\| \nab X^n_s \|_2^2 )^{2\lm}}\gamma
^{z,j}d s \\
\le\| \Gm\|
\int^t_0{ \|\D X^n_s \|_2^2 \over(1+\| \nab X^n_s \|_2^2 )^{2\lm}}\,d s.
\]
Thus, by Doob's inequality,~(\ref{D<I}) and~(\ref{EI}),
\begin{longlist}
\item[(4)]
\begin{eqnarray*}
\\[-28pt]
E \Bigl[\sup_{0 \le t \le T}|N^n_t|^2 \Bigr] &\le& 4E\lan N^n \ran_T
\le C_T
<\infty.
\end{eqnarray*}
\end{longlist}
We conclude (1) from (2)--(4).
\end{pf}

The following estimate plays a key role in the proof of~(\ref{IJtight})
for $p>2$.
%
\begin{lemma}\label{Lemp>2}
Let
%
%
\begin{eqnarray}\label{p1}
p&>&1,\qquad 2<p_1<\infty\qquad\mbox{if $d = 2$,} \nonumber\\[-8pt]\\[-8pt]
p&>&{3d-4 \over d},\qquad 2< p_1 < p{d \over d-2} \qquad
\mbox{if $d \ge
3$}\nonumber
\end{eqnarray}
and let
%
%
\begin{equation}\label{p2}
p_2 <p/\tht_1\qquad
\mbox{where }
\tht_1 ={ {1 /2} -{1 / p_1}\over{1 / 2} -{(d-2) / (dp)}}
\in(0,1).\vadjust{\goodbreak}
\end{equation}
Then, for any
$\del>0$, there are $b,C \in(0,\infty)$ such that for $v \in\cV$
%
%
\begin{equation}\label{p>2}
\| \nab v \|_{p_1}^{p_2}\le C{\cI_p (v) \over(1+\| \nab v \|_2^2
)^\lm}
+C(1+\| \nab v \|_2^2 )^b(1+\| \nab v \|_p)^\del,
\end{equation}
where $\lm$ is defined by~(\ref{lm}). For $d \ge3$, it is possible
to take $\del=0$.
\end{lemma}
%
\begin{pf}
Let $q=2$ for $d \ge3$, $q\in(1,2)$ for $d=2$,
$p_3=p{d \over d-q}>p$. The choice of $p_3$ is made so that
\begin{longlist}
\item[(1)]
\begin{eqnarray}
\nonumber\\[-28pt]
\| \nab v\|_{p_3} \le
C\cI_p (v)^{q /(2p)}(1+\| \nab v \|_p)^{(2-q)/2}\nonumber\\
&&\eqntext{\mbox{cf.~\cite{MNRR96}, page 227, (3.27)}.}
\end{eqnarray}
\end{longlist}
Note also that the choice of $\tht_1$ in~(\ref{p2}) implies that
\begin{longlist}
\item[(2)]
\begin{eqnarray*}
\\[-28pt]
{1 \over p_1}&=&{1 -\tht_1 \over2}+{\tht_1 \over p_3}.
\end{eqnarray*}
\end{longlist}
With $\beta={1 -\tht_1 \over2}+ {\lm q \over2p}\tht_1$
and an arbitrary $\tht_2 \in(0,1)$, we have that
\begin{eqnarray*}
\| \nab v \|_{p_1}
& \st{(2)}{\le} &
\| \nab v \|_2^{1 -\tht_1}\| \nab v \|_{p_3}^{\tht_1}\\
&\st{(1)}{\le}& C\| \nab v \|_2^{1 -\tht_1}
\cI_p (v)^{\tht_1 q /(2p)}(1+\| \nab v \|_p)^{{((2-q)/2)}\tht
_1} \\
& \st{\mathrm{choice}\ \mathrm{of}\ \beta}{=} &
C\biggl( {\cI_p (v) \over(1+\| \nab v \|_2^2)^\lm}\biggr)^{\tht_1 q /(2p)}
(1+\| \nab v \|_2^2)^\beta\\
&&{}\times(1+\| \nab v \|_p)^{{((2-q)/2)}\tht_1} \\
& \st{\tht_2 +(1-\tht_2)=1}{\le} &
C\biggl( {\cI_p (v) \over(1+\| \nab v \|_2^2)^\lm}
\biggr)^{{(\tht_1 /\tht_2)}{(q /(2p))}}\\
&&{}+
C (1+\| \nab v \|_2^2)^{\beta/(1-\tht_2)}
(1+\| \nab v \|_p)^{{((2-q) /2)}{(\tht_1 /(1-\tht_2))}}
\end{eqnarray*}
and hence, that
\[
\| \nab v \|_{p_1}^{p_2}
\le C{\cI_p (v) \over(1+\| \nab v \|_2^2)^\lm}
+C (1+\| \nab v \|_2^2)^b
(1+\| \nab v \|_p)^\del,
\]
where
\[
p_2={\tht_2 \over\tht_1}{2p \over q},\qquad
b={\beta\over1-\tht_2}p_2,\qquad
\del={2-q \over2}{\tht_1 \over1-\tht_2}p_2.
\]
In particular, for $d \ge3$,
\[
p_2={\tht_2 \over\tht_1}p,\qquad
b={\beta\over1-\tht_2}p_2,\qquad
\del=0.
\]
Choosing $\tht_2$ close to 1 (and then $q$ close 2 if $d=2$), we get
the lemma.\vadjust{\goodbreak}
\end{pf}

Lemma~\ref{Lemp>2} is used to obtain the following tightness lemma,
which takes care of the case of $p >2$.
%
\begin{lemma}\label{Lemtightp}
Suppose that
%
%
\begin{eqnarray}\label{p1*}
p &\ge&2,\qquad 2<p_1<\infty\qquad\mbox{if $d = 2$,}
\nonumber\\[-8pt]\\[-8pt]
p &\ge&1+{2d \over d+2},\qquad 2< p_1 < p{d \over d-2} \qquad\mbox{if $d \ge
3$}\nonumber
\end{eqnarray}
and that~(\ref{p2}) holds. Then
%
%
\begin{equation}\label{tightp}
\int^T_0\| \nab X^n_t\|_{p_1}^{p_2}\,dt,\qquad n=1,2,\ldots,\qquad \mbox{are tight}.
\end{equation}
%
\end{lemma}
%
\begin{pf}
By~(\ref{p>2}),
\begin{eqnarray*}
\int^T_0\| \nab X^n_t \|_{p_1}^{p_2}\,dt
& \le&
C\int^T_0{\cI_p (X^n_t) \over(1+\| \nab X^n_t \|_2^2 )^\lm}\,dt \\
& &{} +C\sup_{0 \le t \le T}(1+\| \nab X^n_t \|_2^2 )^b
\int^T_0(1+\| \nab X^n_t \|_p )^\del \,dt.
\end{eqnarray*}
The random variables on the right-hand side ($n=1,2,\ldots$) are tight,
because of
(\ref{apriori}),~(\ref{EI}) and~(\ref{tight2}).
\end{pf}
\begin{pf*}{Proof of Lemma~\ref{LemS^mn->0}}
As explained earlier
[(\ref{xi+Gm->0}), Lemma~\ref{LemGcpt}],
it is enough to show~(\ref{IJtight}).
We recall from~(\ref{IJ}) that
\begin{eqnarray*}
I^{m,n}_T&=&\int^T_0 \| X^m_s \|_{p,\alpha}\| X^n_s \|_{p,1}\| X^n_s
\|_2\,ds,\\
J^{m,n}_T&=&\int^T_0 \| X^m_s \|_{p,\alpha}(1+\| \nab X^n_s \|_p)^{p-1}\,ds.
\end{eqnarray*}

\textit{Case} $1$ ($p=2$).
\begin{eqnarray*}
I^{m,n}_T & \le& \sup_{0 \le t \le T}\| X^n_t \|_{2,1}^2\int^T_0\|
X^m_t \|_{2,\alpha} \,dt, \\
J^{m,n}_T & \le& \sup_{0 \le t \le T}(1+\| \nab X^n_t \|_2)\int
^T_0\| X^m_t \|_{2,\alpha} \,dt.
\end{eqnarray*}
By~(\ref{317}) and~(\ref{tight2}),
the random variables on the right-hand side ($m,n=1,2,\ldots$) are tight
for some $\alpha>1$.

\textit{Case} $2$ ($2<p < {2d \over d-2}$). As for $I_{m,n}$,
\[
I^{m,n}_T \le\sup_{0 \le t \le T}\| X^n_t \|_2
\biggl( \int^T_0\| X^m_t \|_{p,\alpha}^{p'} \,dt\biggr)^{1/p'}
\biggl( \int^T_0\| \nab X^n_t \|_{p,1}^{p} \,dt\biggr)^{1/p}.
\]
Note that $p' <2<p$, since $p>2$. Thus, by~(\ref{apriori}) and~(\ref{317}),
the random variables on the right-hand side ($m,n=1,2,\ldots$)
are tight for some $\alpha>1$.
As for $J^{m,n}_T$, we take $\tl{p} \in(1,p)$ so close to $p$ that
\[
p_2\st{\mathrm{def}}{=}(p-1){\tl{p} \over\tl{p}-1}<p/\tht_1,
\qquad\mbox{where }
\tht_1 ={ {1 /2} -{1 / p}\over{1 / 2} -{(d-2) / (dp)}}
\in(0,1).
\]
Then
\[
J^{m,n}_T \le
\biggl( \int^T_0\| X^m_t \|_{p,\alpha}^{\tl{p}} \,dt\biggr)^{1/\tl{p}}
\biggl( \int^T_0(1+\| \nab X^n_t \|_p)^{p_2} \,dt\biggr)^{(\tl{p}-1) /\tl{p}}.
\]
By~(\ref{317}) and~(\ref{tightp}),
the random variables on the right-hand side ($m,n=1,2,\ldots$)
are tight for some $\alpha>1$.
\end{pf*}
%
\begin{pf*}{Proof of~(\ref{supZ->0})}
Since $p \ge1+{d \over2}$, or equivalently, ${2p \over2p-d } \le p$,
%
%
\begin{equation}\label{R^mTtight}
R^m_T\st{\mathrm{def}}{=}\int^T_0 \|\nab X^m_s\|_p^{2p /(2p-d)}\,ds,
\qquad m=1,2,\ldots,\qquad
\mbox{are tight}
\end{equation}
by~(\ref{apriori}),
and so are $\exp(CR^m_T)$, $m=1,2,\ldots.$ Thus, by Lemma~\ref{LemS^mn->0},
\[
\sup_{0 \le t \le T}\|Z^{m,n}_t\|_2^2 \le S^{m,n}_T\exp(CR^m_T)
\st{m,n \longrightarrow\infty}{\longrightarrow} 0 \qquad\mbox{in
probability}.
\]
Therefore, we get~(\ref{supZ->0}) for $\alpha=0$.
We get~(\ref{supZ->0}) for $\alpha\in(0,1)$ by
interpolation and~(\ref{tight2}).
\end{pf*}

\subsection{\texorpdfstring{Proof of Corollary \protect\ref{Corwsol}}{Proof of Corollary 2.2}} \label{pCor}
We only have to prove~(\ref{SNS}) and~(\ref{apriori2}).
By~(\ref{sdeGar^n}) and integration by parts,
we have for all $\vp\in\cV$ and $t \ge0$
%
%
\begin{eqnarray} \label{SNSn}
\lan\vp, X^n_t \ran
&=& \lan\cP_n\vp, \xi\ran
+\int^t_0 \bigl(\lan X^n_s, (X^n_s \cdot\nab)\vp\ran
- \lan e (\vp), \tau(X^n_s) \ran\bigr)\,ds\nonumber\\[-8pt]\\[-8pt]
&&{}+\lan\cP_n\vp, W_t \ran.\nonumber
\end{eqnarray}
Now, we have by~(\ref{supZ->08}) that
\[
\sup_{0 \le t \le T}|\lan\vp, X^n_t-X_t \ran|\st{n \nearrow\infty
}{\longrightarrow}0\qquad
\mbox{in probability}.
\]
On the other hand, we have by~(\ref{intp1->08}) and the
argument of~\cite{TY10}, Lemma 4.1.1, that
\begin{eqnarray*}
\int^T_0 |\lan X^n_t, (X^n_t \cdot\nab)\vp\ran
-\lan X_t, (X_t \cdot\nab)\vp\ran| \,dt
&\st{n \nearrow\infty}{\longrightarrow}& 0
\qquad\mbox{in probability}, \\
\int^T_0 |\lan e( \vp), \tau(X^n_t)-\tau( X_t) \ran|\,dt
&\st{n \nearrow\infty}{\longrightarrow}&0
\qquad\mbox{in $L_1(P)$}.
\end{eqnarray*}
Therefore, we get~(\ref{SNS}) via~(\ref{SNSn}).
The bound~(\ref{apriori2}) follows from~(\ref{apriori}) by
Fatou's lemma.

\section{\texorpdfstring{Proof of Theorem \protect\ref{Thmee}}{Proof of Theorem 2.3}}
\label{pee}
\subsection{The strategy}
Note that
\[
\| X^n_t\|_2^2 =\sum_{z,j} |X^{n,z,j}|^2.
\]
Applying It\^o's formula to $|X^{n,z,j}|^2$
and using~(\ref{sdeGar^n2}), we see that
%
%
\begin{eqnarray}\label{|X|^2}
|X^{n,z,j}_t|^2&=&|X^{n,z,j}_0|^2+2\int^t_0X^{n,z,j}_s\,dW^{z,j}_s
+2\int^t_0X^{n,z,j}_sb^{z,j}_s(X^n_s)\,ds\nonumber\\[-8pt]\\[-8pt]
&&{}+\lan\psi_{z,j},\Gm\psi_{z,j}\ran t.\nonumber
\end{eqnarray}
Thus,
\[
\| X^n_t \|_2^2-\| X^n_0 \|_2^2= 2M^n_t
+2\int^t_0 \lan X^n_s, b(X^n_s) \ran \,ds +
\operatorname{tr} (\Gm\cP_n)t,
\]
where
%
%
\begin{equation}\label{M^n}
M^n_t=\sum_{z,j}\int^t_0 X^{n,z,j}_s \,dW^{z,j}_s.
\end{equation}
We now recall that
%
%
\begin{equation}\label{uu=0}
\lan w, (v \cdot\nab)w \ran=0,
\end{equation}
$v \in\cV$ and $w \in C^1 (\T^d \ra\rd)$.
Since
\[
\lan v, b (v)\ran
\st{\mbox{{\fontsize{8.36pt}{10.36pt}\selectfont{(\ref{bv2}),
(\ref{uu=0})}}}}{=}
- \lan\tau(v), e (v)\ran,
\]
we have
%
%
\begin{equation}\label{X^2Ito}
\| X^n_t \|_2^2-\| X^n_0 \|_2^2= 2M^n_t
-2\int^t_0 \lan\tau( X^n_s ), e(X^n_s) \ran \,ds +
\operatorname{tr} (\Gm
\cP_n)t.
\end{equation}
Thus, Theorem~\ref{Thmee} follows from the following two lemmas.
%
\begin{lemma}\label{LemconvM}
Referring to~(\ref{M^n}), there exists a martingale $M$ such that
%
%
\begin{equation}\label{convM}
\lim_{n \ra\infty} E \Bigl[ \sup_{0 \le t \le T}|M^n_t-M_t|\Bigr]=0
\qquad\mbox{for any $T \in(0,\infty)$.}
\end{equation}
\end{lemma}
%
%
\begin{lemma}\label{LemconvInt}
For any $T \in(0,\infty)$,
%
%
\begin{equation}\label{convInt}
\int^T_0 | \lan e(X^n_s),\tau(X^n_s) \ran-\lan e(X_s),\tau(X_s)
\ran|\,ds
\st{n \nearrow\infty}{\longrightarrow} 0\qquad
\mbox{in probability $(P)$.}\hspace*{-35pt}
\end{equation}
%
\end{lemma}
%
\subsection{\texorpdfstring{Proof of Lemma \protect\ref{LemconvM}}{Proof of Lemma 4.1}}
It is enough to show that
\begin{longlist}
\item[(1)]
\begin{eqnarray*}
\\[-28pt]
\lim_{m,n \ra\infty} E \Bigl[ \sup_{0 \le t \le T}|M^n_t-M^m_t|\Bigr]&=&0.
\end{eqnarray*}
\end{longlist}
By the Burkholder--Davis--Gundy inequality,
\[
E \Bigl[ \sup_{0 \le t \le T}|M^n_t-M^m_t|\Bigr]
\le C E [ \lan M^n-M^m \ran_T^{1/2} ].
\]
We may assume $m > n$. Then, for any $t>0$,
\[
M^m_t-M^n_t=\mathop{\sum_{z,j}}_{n< |z|_\infty\le m}\int^t_0 X^{m,z,j}_s
\,dW^{z,j}_s
+\mathop{\sum_{z,j}}_{|z|_\infty\le n}\int^t_0 (X^{m,z,j}_s-X^{n,z,j}_s
)\,dW^{z,j}_s
\]
and thus,
\begin{eqnarray*}
\lan M^m-M^n \ran_t
&=&
\mathop{\sum_{z,j}}_{n< |z|_\infty\le m}\int^t_0( X^{m,z,j}_s )^2\gamma^{z,j}\,ds
\\
&&{}+\mathop{\sum_{z,j}}_{|z|_\infty\le n}\int^t_0( X^{m,z,j}_s-X^{n,z,j}_s
)^2\gamma^{z,j}\,ds \\
&\le& Q_t +R_t,
\end{eqnarray*}
where
\[
Q_t=\int^t_0 \bigl\|(1 -\cP_n)\sqrt{\Gm}X^m_s\bigr\|_2^2\,ds,\qquad
R_t=\int^t_0 \bigl\|\sqrt{\Gm}(X^m_s-X^n_s)\bigr\|_2^2\,ds.
\]
By~(\ref{apriori}), we have
\[
E [ Q_T]
\le\bigl\|(1 -\cP_n)\sqrt{\Gm} \bigr\|_{2 \ra2}^2
\int^T_0 E [ \|X^m_s\|
_2^2 ]\,ds
\le\bigl\|(1 -\cP_n)\sqrt{\Gm} \bigr\|_{2 \ra2}^2C_T
\st{ n \nearrow\infty}{\longrightarrow} 0.
\]
On the other hand, we see from~(\ref{supZ->08}) that
\[
E [ R_T^{1/2}]
\le\bigl\|\sqrt{\Gm}\bigr\|_{2 \ra2}
E \biggl[ \biggl( \int^T_0\|X^m_s -X^n_s\|_2^2 \,ds \biggr)^{1/2}\biggr]
\st{ m,n \nearrow\infty}{\longrightarrow} 0.
\]
Putting things together, we get (1).

\subsection{\texorpdfstring{Proof of Lemma \protect\ref{LemconvInt}}{Proof of Lemma 4.2}}
We write
\begin{eqnarray*}
&&\lan\tau(X_s), e(X_s) \ran-\lan\tau(X^n_s), e(X^n_s)
\ran\\
&&\qquad= \lan\tau(X_s)-\tau(X^n_s), e(X_s) \ran
+\lan\tau(X^n_s), e(X_s)-e(X^n_s) \ran.
\end{eqnarray*}
In view of this, we will prove~(\ref{convInt}) by showing that
\begin{longlist}
\item[(1)]
\begin{eqnarray*}
\\[-28pt]
 \int^T_0 |\lan\tau(X_s)-\tau(X^n_s), e(X_s) \ran|\,ds &\st{n
\nearrow\infty}{\longrightarrow}& 0\qquad
\mbox{in probability ($P$),}
\end{eqnarray*}
\item[(2)]
\begin{eqnarray*}
\\[-28pt]
\int^T_0 |\lan\tau(X^n_s), e(X_s)-e(X^n_s)\ran|\,ds &\st{n
\nearrow\infty}{\longrightarrow}& 0
\qquad\mbox{in probability ($P$).}
\end{eqnarray*}
\end{longlist}
To show (1), we note that
\begin{eqnarray}
\bigl| (1+|x|^2)^{(p-2)/2}x-(1+|y|^2)^{(p-2)/2}y \bigr|
\le C|x-y|(1+|x|+|y|)^{p-2},\nonumber\\
&&\eqntext{x,y \in\rd\otimes\rd.}
\end{eqnarray}
Therefore, with $p_1 \in(1,p)$ and $p_1'={p_1 \over p_1-1}$,
\begin{eqnarray*}
&&\int^T_0 |\lan\tau(X_s)-\tau(X^n_s), e(X_s) \ran|\,ds \\
&&\qquad \le
C\int^T_0ds \int_{\T^d}|e (X_s) - e (X^n_s )| \bigl(1+|e (X^n_s )|+|e
(X_s )|\bigr)^{p-1}\\
&&\qquad\le CI_n^{1 / p_1}(I_n')^{1 / p_1'},
\end{eqnarray*}
where
\begin{eqnarray*}
I_n&=& \int^T_0ds \int_{\T^d} | e (X_s) -e (X^n_s)|^{p_1},\\
I'_n&=& \int^T_0 ds \int_{\T^d}\bigl(1+| e (X^n_s)|+| e ( X_s )|\bigr)^{(p-1)p_1'}.
\end{eqnarray*}
Note that $(p-1)p_1' \searrow p$ as $p_1 \nearrow p$. Thus,
for $p_1$ sufficiently close to $p$, $\{I_n'\}_{n \ge1}$
are tight by~(\ref{tightp}).
On the other hand,
$I_n \ra0$ in probability ($P$) for any $p_1 <p$
by~(\ref{intp1->08}). Thus, we get (1).

As for (2), with $p_1 \in(1,p)$ and $p_1'={p_1 \over p_1-1}$ again,
\[
\int^T_0 ds \int_{\T^d}| \lan\tau(X^n_s), e(X_s)-e(X^n_s) \ran|\,ds
\le J_n^{1/p_1}(J_n')^{1/p_1'},
\]
where
\begin{eqnarray*}
J_n &=& \int^T_0ds\int_{\T^d} | e(X_s)-e(X^n_s) |^{p_1}, \\
J'_n &=& \int^T_0 ds\int_{\T^d} | \tau(X^n_s)|^{p_1'}
\le\nu\int^T_0 ds\int_{\T^d}\bigl(1+| e (X^n_s)|\bigr)^{(p-1)p_1'}.
\end{eqnarray*}
As in the proof of (1), for $p_1$ sufficiently close to $p$,
$\{J_n'\}_{n \ge1}$ are tight by~(\ref{tightp}),
and $J_n \ra0$ in probability ($P$) for any $p_1 <p$
by~(\ref{intp1->08}). Thus, we get (2).

\section{\texorpdfstring{Proof of Proposition \protect\ref{Prop2DNS}}{Proof of Proposition 2.4}} \label{2DNS}
\subsection{\texorpdfstring{Proofs of (\protect\ref{2Dtight}) and (\protect\ref{2Dconv})}
{Proofs of (2.23) and (2.24)}}
Note that for $\alpha=0,1,2,\ldots$
\[
\| \nab^\alpha v\|_2^2 =\lan v, (-\D)^\alpha v \ran
=\sum_{z,j} (-4\pi^2 |z|^2)^\alpha\lan v, \psi_{z,j} \ran^2,\qquad
v \in\cV.
\]
By plugging $v=X^n_t$ into the above identity,
and using~(\ref{|X|^2}), we obtain that
%
%
\begin{eqnarray}\label{Ito^2}
\| \nab^\alpha X^n_t \|_2^2&=&\| \nab^\alpha X^n_0 \|_2^2+ 2M^n_t
+2\int^t_0 \lan(-\D)^\alpha X^n_s, b(X^n_s) \ran
\,ds\nonumber\\[-8pt]\\[-8pt]
&&{}+\operatorname{tr} (\Gm(-\D)^\alpha\cP_n)t,\nonumber
\end{eqnarray}
where
%
%
\begin{equation}\label{Ito^2M}
M^n_t = \sum_{z,j}\int^t_0 (-\D)^\alpha X^{n,z,j}_s \,d W^{z,j}_s.
\end{equation}
Since we assume~(\ref{trGD^a}),
we may repeat the proof of Lemma~\ref{LemconvM}, with $\Gm$ replaced by
$\Gm(-\D)^\alpha$ to obtain the following lemma.
%
\begin{lemma}\label{LemconvM2}
Referring to~(\ref{Ito^2M}), there exists a martingale $M$ such that
%
%
\begin{equation}\label{convM2}
\lim_{n \ra\infty} E \Bigl[ \sup_{0 \le t \le
T}|M^n_t-M_t|\Bigr]=0\qquad
\mbox{for any $T \in(0,\infty)$.}
\end{equation}
%
\end{lemma}

We now continue on~(\ref{Ito^2}).
For $p=2$, we have for $v \in\cV$ that
%
%
\begin{eqnarray} \label{b2}
\lan(-\D)^\alpha v, b(v) \ran
& = & \lan(-\D)^\alpha v, (v \cdot\nab) v\ran
+\nu\lan(-\D)^\alpha v, \D v \ran\\
& = & \lan(-\D)^\alpha v, (v \cdot\nab) v\ran
-\nu\| \nab^{\alpha+1} v \|_2^2.\nonumber
\end{eqnarray}
Moreover, we have for $d=2$ that
%
%
\begin{equation}\label{227}
|\lan(-\D)^\alpha v, (v \cdot\nab) v\ran|
\le C_1 \| \nab^{\alpha+1}v \|_2^{(2\alpha-1) /\alpha}\| \nab v
\|_2^{(\alpha+1) /\alpha}.
\end{equation}
This follows from the argument in the proof of (2.27),
\cite{Kuk06}, page 17. By~(\ref{227}) and Young inequality,
we obtain that
%
%
\begin{equation}\label{2272}\qquad
|\lan(-\D)^\alpha v, (v \cdot\nab) v\ran|
\st{{(2\alpha-1)/(2\alpha)}+{1 / (2\alpha)}=1}{\le}
{\nu\over2}\| \nab^{\alpha+1}v \|_2^2+C_2\| \nab v \|_2^{2\alpha+2}.
\end{equation}
By~(\ref{Ito^2}) and~(\ref{b2}),
%
%
\begin{eqnarray}
\label{Ito^22}
& & \| \nab^\alpha X^n_t \|_2^2
+2\nu\int^t_0 \| \nab^{\alpha+1}X^n_s \|_2^2\,ds \nonumber\\
& &\qquad = \| \nab^\alpha X^n_0 \|_2^2+ 2M^n_t +2\int^t_0 \lan(-\D)^\alpha X^n_s, (X^n_s \cdot\nab)
X^n_s\ran \,ds\\
& &\qquad\quad{}
+\operatorname{tr} (\Gm(-\D)^\alpha\cP_n)t. \nonumber
\end{eqnarray}
Therefore, by~(\ref{2272}),
%
%
\begin{eqnarray}\label{Ito^23}
&&\| \nab^\alpha X^n_t \|_2^2
+\nu\int^t_0 \| \nab^{\alpha+1}X^n_s \|_2^2\,ds\nonumber\\
&&\qquad\le \| \nab^\alpha X^n_0 \|_2^2+ 2M^n_t
+C_2 \int^t_0\| \nab X^n_s \|_2^{2\alpha+2}
\\
&&\qquad\quad{}
+\operatorname{tr} (\Gm(-\D)^\alpha\cP_n)t.\nonumber
\end{eqnarray}
We conclude the tightness~(\ref{2Dtight}) from~(\ref{Ito^23}), using
(\ref{trGD^a}),~(\ref{xi2a}), Lemmas~\ref{Lemtight2} and~\ref{LemconvM2}.
The convergence~(\ref{2Dconv}) follows from~(\ref{supZ->08}) and
(\ref{2Dtight}) via interpolation.

\subsection{The pathwise balance relation for the enstrophy}
Here, we prove that the process defined by~(\ref{Mt2}) is a martingale.
Since
\[
\lan\D v, (v \cdot\nab) v\ran=0\qquad
\mbox{for $v \in\cV$}\qquad \mbox{cf.~\cite{MNRR96}, page 225, (3.20)},
\]
we set $\alpha=1$ in~(\ref{Ito^22}) to get
%
%
\begin{equation}\label{Ito^22a=1}\qquad
\| \nab X^n_t \|_2^2
+2\nu\int^t_0 \| \D X^n_s \|_2^2\,ds
= \| \nab X^n_0 \|_2^2+ 2M^n_t
+\operatorname{tr} (\Gm(-\D)\cP_n)t,
\end{equation}
where
\[
M^n_t =\sum_{z,j}\int^t_0 (-\D) X^{n,z,j}_s \,d W^{z,j}_s,
\]
for which Lemma~\ref{LemconvM2} (with $\alpha=1$) is valid.
Since we assume~(\ref{trGD^a}) and~(\ref{xi2a}) with $\alpha=2$,
we have by~(\ref{2Dconv}) that
%
%
\begin{equation}\label{2Dconv2}
\sup_{0 \le t \le T}\| X^n_t -X_t \|_{2,1}^2
+\int^T_0 \| X^n_t -X_t\|_{2,2}^2\,dt\st{n \nearrow\infty
}{\longrightarrow} 0\qquad
\mbox{in probability}.\hspace*{-35pt}
\end{equation}
Therefore, we let $n \nearrow\infty$ in~(\ref{Ito^22a=1})
to see that
\[
\| \nab X_t \|_2^2
+2\nu\int^t_0 \| \D X_s \|_2^2\,ds
= \| \nab X_0 \|_2^2+ 2M_t
+\operatorname{tr} (\Gm(-\D))t, \qquad t \ge0.
\]
This means that the process $M_\cdot$ defined by~(\ref{Mt2}) is
exactly the
martingale obtained in Lemma~\ref{LemconvM2} (with $\alpha=1$).

\section*{Acknowledgments}

The author thanks Professors Franco Flandoli, Reika Fukuizumi and Kenji
Nakanishi for useful conversation.


%

\printaddresses

\end{document}